\documentclass[12pt, reqno]{amsart}
\usepackage{amssymb,amsthm,amsfonts,amsmath, amscd}

\usepackage[hypertex]{hyperref}
\usepackage{mathrsfs}

\newcommand\R{\mathbb R}
\newcommand\Z{\mathbb Z}
\newcommand\C{\mathbb C}
\newcommand\T{\mathbb T}

\newcommand\al{\alpha}
\newcommand\be{\beta}
\newcommand\ga{\gamma}

\newcommand\Ga{\Gamma}
\newcommand\de{\delta}
\newcommand\Om{\Omega}
\newcommand\om{\omega}

\newcommand\la{\lambda}
\newcommand\La{\Lambda}
\newcommand\epsi{\varepsilon}

\newcommand\B{\mathscr B}
\newcommand\M{\mathscr M}

\newcommand\X{\mathfrak X}
\newcommand\Y{\mathfrak Y}

\newcommand\const{\operatorname{const}}
\newcommand\tr{\operatorname{tr}}
\newcommand\Dom{\operatorname{Dom}}

\newcommand\pd{\partial}
\newcommand\wt{\widetilde}
\newcommand\wh{\widehat}
\newcommand\zw{{(z,z',w,w')}}
\newcommand\uv{{(u,u',v,v')}}

\newtheorem{theorem}{Theorem}[section]

\theoremstyle{definition}

\newtheorem{remark}[theorem]{Remark}
\newtheorem{example}[theorem]{Example}

\numberwithin{equation}{section}

\begin{document}

\title[]{Markov dynamics on the dual object to the infinite-dimensional unitary group}

\author{Grigori Olshanski}
\address{Institute for Information Transmission Problems, Moscow, Russia;
\newline\indent Independent University of Moscow, Russia;
\newline \indent National Research University Higher School of
Economics, Moscow, Russia} \email{olsh2007@gmail.com}

\thanks{Partially supported by a grant from Simons Foundation
(Simons--IUM Fellowship) and the project SFB 701 of Bielefeld University.}


\maketitle

\tableofcontents

\section{Preface}\label{sect1}

These are notes for a mini-course of 3 lectures given at the St. Petersburg
School in Probability and Statistical Physics (June 2012). My aim was to
explain, on the example of a particular model, how ideas from the
representation theory of big groups can be applied in probabilistic problems.
The material is based on the joint paper \cite{BO-GT-Dyn} by Alexei Borodin and
myself; a broader range of topics is surveyed in the lecture notes by Alexei
Borodin and Vadim Gorin \cite{BG-Lectures}.

The main result of \cite{BO-GT-Dyn} consisted in constructing a family of
Feller Markov processes living on the infinite-dimensional locally compact
space $\wh{U(\infty)}$, a kind of dual object to the infinite-dimensional
unitary group $U(\infty)$. By definition, the group $U(\infty)$ is the union of
the chain of compact unitary groups $U(N)$, $N=1,2,\dots$, embedded to each
other. Dually, the space $\wh{U(\infty)}$ appears as the ``entrance boundary''
of a chain of discrete sets $\wh{U(N)}$ related to each other by certain
stochastic matrices. This structure plays a key role in our construction of
Markov dynamics on $\wh{U(\infty)}$.

The problem solved in \cite{BO-GT-Dyn} is in some (nonconventional and not
strictly defined) sense dual to the problem of constructing an
infinite-dimensional analog of the fundamental Dyson's model \cite{Dys} of an
$N$-particle non-colliding process coming from the Brownian motion on $U(N)$.
The latter problem, initiated by Spohn \cite{Spo87}, is investigated in recent
works Katori--Tanemura \cite{KatoriTanemura-JSP09},
\cite{KatoriTanemura-CMP10}, \cite{KatoriTanemura-MPRF11}, and Osada
\cite{Osada-AnnProb13}. In our problem, the role of Dyson's model is played by
a family of continuous time Markov chains on $\wh{U(N)}$. At first glance, it
looks much more sophisticated than Dyson's model but actually it turns out to
be more friendly.

The method used in \cite{BO-GT-Dyn} was also applied to other models in
Borodin--Gorin \cite{BG-PTRF} and Borodin--Olshanski \cite{BO-ThomaCone}.

The prerequisites for reading the present notes are modest: an acquaintance
with the basics of Markov processes is enough, and no real knowledge of
representation theory is assumed.

I would like to thank Alexei Borodin for valuables comments. I am also grateful
to Stanislav Smirnov for the opportunity to take part in the program of the
School.

\section{Dyson's model}\label{sect2}

Let us start with recalling a classical fact. Consider the classical
multidimensional Brownian motion in $\R^N$, $B\!M(\R^N)$, whose generator is
half the Laplacian. Because the Laplacian admits a separation of variables in
the polar coordinates, the radial part of $B\!M(\R^N)$ is still a Markov
process. Namely, it is the Bessel process $B\!E\!S^N$ on the halfline
$\R_+=\{r\in\R: r\ge0\}$; the generator of $B\!E\!S^N$ is the ordinary
differential operator
$$\
\frac12\left(\frac{d^2}{dr^2}+\frac{N-1}r\frac{d}{dr}\right).
$$
See, e.g., It\^o-McKean \cite{ItoMcKean}.

A similar effect holds for a number of other multidimensional diffusion
processes, in particular, for the {\it Brownian motion on the unitary group\/},
see, e.g., Dyson \cite{Dys}, McKean\cite{McKean}. This diffusion process, which
we denote by $B\!M(U(N))$, lives on the group $U(N)$ of $N\times N$ unitary
matrices and is generated by a two-sided invariant second order differential
operator on that group.  The analog of the radial projection $\R^N\to\R_+$ is
the map assigning to a generic unitary matrix $g\in U(N)$ the collection
$(u_1,\dots,u_N)$ of its eigenvalues, which we interpret as an unordered
$N$-tuple of points on the unit circle $\T:=\{u\in\C: |u|=1\}$. Note that if
$g\in U(N)$ is in general position, then the eigenvalues $u_i$ are pairwise
distinct. The assignment $g\mapsto (u_1,\dots,u_N)$ maps $U(N)$ onto
$\T^N/S_N$, the quotient of the $N$-fold product space $\T^N$ with respect to
the action of the symmetric group $S_N$ permuting the coordinates. Thus,
$\T^N/S_N$ plays the role of the halfline.

It turns out that one can define the radial part of $B\!M(U(N))$, which is a
diffusion process on $\T^N/S_N$; let us denote it by $X_N$.

To describe the generator of $X_N$, it is convenient to pass from the
``multiplicative coordinates'' $u_1,\dots,u_N$ to the ``additive coordinates''
$x_1,\dots,x_N$ by setting $u_k=\exp(\sqrt{-1}\,x_k)$, where $k=1,\dots,N$ and
$x_k\in\R/2\pi\Z$. In these coordinates, the generator in question, denoted by
$D_N$, can be written in the form
\begin{equation}\label{eq2.A}
D_N=V_N^{-1}\circ\Delta_N\circ V_N+C_N,
\end{equation}
where
\begin{gather*}
V_N=V_N(x_1,\dots,x_N):=\prod_{1\le i<j}|u_i-u_j|=\const\prod_{1\le i<j\le
N}\sin\frac{x_i-x_j}2,\\
\Delta_N:=\sum_{i=1}^N\frac{\pd^2}{\pd x_i^2},\\
C_N=\frac{(N-1)N(N+1)}{12}\,
\end{gather*}
In words, \eqref{eq2.A} means that to apply $D_N$ to a function $F$ we first
multiply $F$ by $V_N$, then apply the Laplacian $\Delta_N$, then divide by
$V_N$, and finally add $C_NF$. Since
$$
\Delta_N V_N=-C_NV_N,
$$
$D_N$ annihilates the constants.

Formula \eqref{eq2.A} is a kind of {\it Doob's $h$-transform\/} (see
Rogers--Williams \cite{RogersWilliams}) applied to the ``flat'' Brownian motion
generated by the Laplacian $\Delta_N$, where $h=V_N$. However, $V_N$ is not a
harmonic function for $\Delta_N$ but only an eigenfunction; this explains the
appearance of the compensating term $C_N$.

More explicitly, \eqref{eq2.A} can be rewritten as
\begin{equation}\label{eq2.B}
D_N=\sum_{i=1}^N\frac{\pd^2}{\pd x_i^2}+\sum_{i=1}^N\left(\sum_{\al:\, \al\ne
i}\cot\frac{x_i-x_\al}{2}\right)\frac{\pd}{\pd x_i}.
\end{equation}
Although the coefficients of the first order derivatives have singularities
along the diagonals $x_i=x_j$, these singularities are cancelled when $D_N$ is
applied to smooth {\it symmetric\/} functions in $x_1,\dots,x_N$. Note that
natural ``observables'' on the quotient space $\T^N/S_N$ are just symmetric
functions in the coordinates on $\T^N$.

In contrast to the Bessel process, the process $X_N$ generated by $D_N$ has a
stationary distribution $\mu_N$; it is the radial part of the normalized Haar
measure on $U(N)$. The density of $\mu_N$ with respect to Lebesgue measure
$dx_1\dots dx_N$ is proportional to $V_N^2$. The measure $\mu_N$ first emerged
in the context of Weyl's character formula, see Weyl \cite{Wey39}. In Dyson's
interpretation, $\mu_N$ is the law of a system of $N$ interacting point
particles on the unit circle $\T^N$, called the  {\it circular unitary
ensemble\/} and usually denoted as $C\!U\!E_N$, see Dyson \cite{Dys}, Mehta
\cite{Mehta}.

A well-known result says that, in a suitable large-$N$ scaling limit regime,
$C\!U\!E_N$ turns into an ensemble of infinitely many interacting particles on
$\R$. The distribution of the particles is best described in terms of the
correlation functions, which are determined by a simple translation invariant
correlation kernel on $\R$, called the sine kernel. See, e.g., Mehta
\cite{Mehta}.

In view of this fact it is natural to ask what happens with the process $X_N$
in the same scaling limit regime: does there exists a Markov process $X_\infty$
which would be a large-$N$ limit (in some reasonable sense) of the diffusions
$X_N$? Using the dynamical correlation functions one can check that if the
initial distribution of the process $X_N$ is $\mu_N$, then its multi-time
finite-dimensional distributions survive in a suitable scaling limit transition
as $N\to\infty$. However, this is insufficient to conclude that $X_\infty$ does
exist.

In the attempt to imagine the possible form of the generator of $X_\infty$, let
us examine the limiting behavior of the operators $D_N$ (the informal argument
below follows the discussion in the beginning of Spohn's paper \cite{Spo87}).

The scaling limit in question consists in a change of variables,
$$
x_i\leadsto y_i, \qquad y_i:=\frac{N}{2\pi}\,x_i, \qquad 1\le i\le N.
$$
We assume that the initial coordinates $x_i$ range over the interval
$(-\pi,\pi)$; then the new coordinates $y_i$ range over the interval
$(-N/2,N/2)$ of length $N$, so that the mean density of particles equals 1.
Writing
$$
\cot\frac{x_i-x_\al}2=\cot\left(\frac\pi N(y_i-y_\al)\right)\sim\frac N{2\pi}
\frac2{y_i-y_\al}
$$
we get
$$
\left(\frac{2\pi}{N}\right)^2 D_N \sim\sum_i\frac{\pd^2}{\pd
y_i^2}+2\left(\sum_{\al:\, \al\ne i}\frac1{y_i-y_\al}\right)\frac{\pd}{\pd
y_i}.
$$
Next, the factor $\left(\frac{2\pi}{N}\right)^2$ can be eliminated by rescaling
the time parameter, so that finally we are left with the formal differential
operator
$$
D_\infty:=\sum_i\frac{\pd^2}{\pd y_i^2}+2\left(\sum_{\al:\, \al\ne
i}\frac1{y_i-y_\al}\right)\frac{\pd}{\pd y_i},
$$
where we may assume that the index $i$ ranges over $\Z$ and
$\ldots<y_{-1}<y_0<y_1<\ldots$.

We see that the series in the brackets, in general, diverges, and even if we
manage to regularize it, it is highly non-evident how to prove that a suitable
regularization of $D_\infty$ serves as a (pre)generator of a Markov process.

This informal argument shows that a naive direct approach to constructing the
generator of $X_\infty$ faces serious difficulties.

\section{The one-particle dynamics: a bilateral birth-death
process}\label{sect3}

I proceed to the  model from our work \cite{BO-GT-Dyn}. It bears some
resemblance with Dyson's model and (in some informal sense to be clarified
below) is dual to it.

Let us start with the simplest case when $N$, the number of particles, equals
1. The one-particle Dyson model is very simple, it is the conventional Brownian
motion on the unit circle $\T$. I will explain what I mean by the corresponding
``dual model''.

The space $\T$ is a compact Abelian group, and its Pontryagin dual is the
discrete Abelian group $\Z$. So it is not surprising that the ``dual'' state
space is $\Z$. But what is  a substitute of the Brownian motion? As $\Z$ is
discrete, it cannot be a diffusion process, it should be a jump process or, in
other words, a Markov chain. We want a continuous time process, so that it is a
continuous time Markov chain.

The generator of the Brownian motion on $\T$ is the simplest second order
differential operator, $d^2/dx^2$, where $x$ is the ``additive coordinate'' as
above. A natural lattice analog of this operator should be a second order
difference operator $D$ on $\Z$ transforming a test function $F(l)$ to the
function
\begin{equation}\label{eq3.A}
DF(l)=a^+(l)(F(l+1)-F(l))+a^-(l)(F(l-1)-F(l)), \qquad l\in\Z,
\end{equation}
where the coefficients $a^+(l)$ and $a^-(l)$ represent the rates of the jumps
$l\to l+1$ and $l\to l-1$, respectively (these are the only possible jumps).

At first glance, the most natural choice of the coefficients is to set
$a^+(l)=a^-(l)=\const>0$. This leads to the Markov chain which looks as the
most natural discrete analog of the classical Brownian motion. However, this
chain is not suitable for our purposes as it does not possess a stationary
distribution. Certainly, the counting measure on $\Z$ is invariant, but it is
infinite, while we would like to have a finite measure, as in the case of $\T$.
For this reason we reject the constant coefficients.

The next possible variant would be  to make the coefficients $a^\pm(l)$ some
linear functions in $l$. This indeed allows one to get examples of processes
with a stationary distribution (some birth-death processes). However, they
cannot live on the whole lattice $\Z$, because a linear function changes the
sign, while the jump rate cannot take negative values. Since we want to deal
with the whole lattice, we reject this variant as well.

Let us try now quadratic rates $a^\pm(l)$. The leading terms in $a^+(l)$ and
$a^-(l)$ must coincide to prevent a growing drift to $+\infty$ or $-\infty$
(such a drift is obviously incompatible with a stationary distribution). Then,
without loss of generality,  we may assume that $a^\pm(l)$ equals $l^2$ plus
lower degree terms. Writing such a quadratic function as the product of two
linear factors we set
\begin{equation}\label{eq3.B}
a^+(l)=(u-l)(u'-l), \quad a^-(l)=(v+l)(v'+l),
\end{equation}
where $(u,u',v,v')$ is a quadruple of parameters. Note that the change $l\to
-l$ amounts to switching $(u,u')\leftrightarrow(v,v')$.

Finally, we want $a^\pm(l)$ to take strictly positive values for all $l\in\Z$.
Let us say that a couple $(z,z')$ of complex numbers is {\it admissible\/} if
$(u-l)(u'-l)>0$ for all $l\in\Z$. We will assume that both $(u,u')$ and
$(v,v')$ are admissible.

It is not difficult to classify all admissible couples. Namely, $(u,u')$ is
admissible if and only if

\medskip

$\bullet$ either both $u$ and $u'$ are nonreal complex numbers and $u'=\bar u$;

$\bullet$ or both $u$ and $u'$ are real and there exists $m\in\Z$ such that
$m<u,u'<m+1$.

\medskip

It turns out that quadratic rates give the desired result:

\begin{theorem}\label{thm3.A}
Assume $(u,u')$ and $(v,v')$ are admissible couples of parameters.

{\rm(i)} There exists a continuous time Markov chain $X^\uv_1$ on $\Z$, such
that the only possible  jumps are of the form $l\to l\pm1$ and their rates
$a^\pm(l)$ are given by \eqref{eq3.B}.

{\rm(ii)} The chain $X^\uv_1$ possesses a unique, within a constant factor,
symmetrizing measure. This measure is finite if and only if  the parameters
satisfy the additional constraint $u+u'+v+v'>-1$.
\end{theorem}

Note that $u+u'+v+v'$ is a real number because $(u,u')$ and $(v,v')$ are
admissible. Note also that a symmetrizing measure is automatically invariant.
We denote the symmetrizing measure of our Markov chain by $M^\uv_1$. Here is an
explicit expression for it:
\begin{equation}\label{eq3.C}
M^\uv_1(l)=\const\,\frac1{\Ga(u+1-l)\Ga(u'+1-l)\Ga(v+1+l)\Ga(v'+1+l)}\,.
\end{equation}

The normalization constant is found from a beautiful classical hypergeometric
identity due to Dougall \cite{Dougall} (see also Erdelyi \cite[\S1.4]{Er}),
\begin{multline}\label{eq3.D}
\sum_{l\in\Z}\frac1{\Ga(u+1-l)\Ga(u'+1-l)\Ga(v+1+l)\Ga(v'+1+l)}\\
=\frac{\Ga(u+u'+v+v'+1)}{\Ga(u+v+1)\Ga(u+v'+1)\Ga(u'+v+1)\Ga(u'+v'+1)}\,.
\end{multline}
If $u+u'+v+v'>-1$, then we may take as the constant factor in \eqref{eq3.C} the
quantity inverse to the right-hand side in \eqref{eq3.D}; with this
normalization  $M^\uv_1$ becomes a probability measure.

The Markov chain  $X^\uv_1$ is an example of so-called bilateral birth and
death processes, see Feller \cite[Section 17]{Fel57}, Pruitt \cite{Pru63}, Yan
\cite{Yan90}.

The above arguments are intended to convince the reader that the definition of
the chain $X^\uv_1$ is quite natural. But in reality, this definition came from
other considerations, related to our previous work on harmonic analysis on big
groups: \cite{BO-AnnMath}, \cite{BO05-EurCongr}, \cite{BO-PTRF}, \cite{Ols03}.

\section{The $N$-particle dynamics}\label{sect4}

In this section $\uv$ is a fixed quadruple of parameters such that $(u,u')$ and
$(v,v')$ are admissible, and $N\ge2$ is a fixed natural number.

We are dealing with the lattice $\Z^N$; its elements are denoted as
$\ell=(l_1,\dots,l_N)$. Denote by $D$ the 1-dimensional difference operator
\eqref{eq3.A} with the coefficients given by \eqref{eq3.B}, and let $D^{[i]}$
stand for a copy of $D$ acting on the $i$th coordinate of $\ell$, where
$i=1,\dots,N$. The operator
\begin{equation}\label{eq4.A}
D_N^{\rm free}:=\sum_{i=1}^N D^{[i]}.
\end{equation}
generates a continuous time Markov chain on $\Z^N$, which is simply the product
of $N$ independent copies of the chain $X_1^\uv$.

The next step is to apply to this chain the Doob $h$-transform (cf.
\eqref{eq2.A}), taking as $h$ the function
$$
V_N(\ell):=\prod_{1\le i<j\le N}(l_i-l_j).
$$
Note that $V_N$ is an eigenfunction of $D_N^{\rm free}$:
$$
D_N^{\rm free}V_N=-C_N V_N,
$$
where
$$
C_N=C_N\uv:=\frac{N(N-1)}2(u+u'+v+v')-\frac{N(N-1)(N-2)}3\,.
$$

Consider the region $\Om_N\subset\Z^N$ defined by
\begin{equation}\label{eq4.C}
\Om_N:=\{\ell\in\Z^N: \,  l_1>\dots>l_N\}.
\end{equation}

\begin{theorem}\label{thm4.A}
Assume $(u,u')$ and $(v,v')$ are admissible couples of parameters.

{\rm(i)} There exists a continuous time Markov chain $\wt X^\uv_N$ on\/
$\Om_N\subset\Z^N$, whose infinitesimal generator is given by the difference
operator
$$
D_N:=V_N^{-1}\circ D_N^{\rm free}\circ V_N+C_N.
$$

{\rm(ii)} The chain $\wt X^\uv_N$ possesses a unique, within a constant factor,
symmetrizing {\rm(}and hence invariant{\rm)} measure $\wt M^\uv_N$:
\begin{multline}\label{eq4.D}
\wt M^\uv_N(\ell)=\const \prod_{i=1}^N
\bigg(\frac1{\Gamma(u-\ell_i+1)\Gamma(u'-\ell_i+1)}\\
\times\frac1{\Gamma(v+\ell_i+1)\Gamma(v'+\ell_i+1)}\bigg)\cdot (V_N(\ell))^2,
\qquad \ell\in\Om_N.
\end{multline}
This measure is finite if and only if  the parameters satisfy the additional
constraint $u+u'+v+v'>2N-3$.
\end{theorem}

The possible jumps of the chain are of the form $\ell\to\ell\pm\epsi_i$, where
$\epsi_i$ denotes the $i$th basis vector in $\R_N$, $i=1,\dots,N$. Note that if
$\ell\in\Om_N$ but $\ell+\epsi_i$ or $\ell-\epsi_i$ does not belong to $\Om_N$,
then the corresponding rate automatically vanishes, so that the chain does not
exit from $\Om_N$.

Explicitly, the rate of the jump $\ell\to\ell\pm\epsi_i$ equals
\begin{equation}\label{eq4.B}
\frac{V_N(\ell\pm\epsi_i)}{V_N(\ell)}\cdot\left\{\begin{matrix}
(u-l_i)(u'-l_i)\\(v+l_i)(v'+l_i)\end{matrix}\right\},
\end{equation}
where the upper/lower quantity in the braces corresponds to the plus/minus
sign, respectively.

As will be shown in Section \ref{sect7}, $\Om_N$ serves as the set of
parameters for $\wh{U(N)}$, the dual object to the unitary group $U(N)$. This
is why we view the dynamics just introduced as ``dual'' to the Dyson model of
Section \ref{sect2}.

\section{The method of intertwiners}\label{sect5}

Here I describe a general formalism which will be used for constructing a model
of infinite-dimensional Markov dynamics out of the Markov chains $\wt X^\uv_N$.
For more detail, see Borodin--Olshanski \cite{BO-GT-Dyn}, \cite{BO-MMJ}.

An $m$-dimensional simplex $\Delta^m$ in a vector space has $m+1$ vertices, and
each point of $\Delta^m$ is uniquely represented as a convex linear combination
of the vertices. It follows that $\Delta^m$ can be identified with the set of
the probability measures on the set of the vertices.

Let us adopt this viewpoint and, more generally, given a finite or countably
infinite abstract set $\X$, we define the simplex with the vertex set $\X$ as
the set of probability measures on $\X$.

Even more generally, let $\X$ be a {\it measurable space\/}, that is, a set
with a distinguished $\sigma$-algebra $\B(\X)$ of subsets called measurable
subsets (in another terminology, $\X$ is a {\it Borel space\/}). We assume that
$\B(\X)$ contains all singletons. Denote by $\M(\X)$ the space of probability
measures defined on $\B(\X)$. We regard $\M(\X)$ as a {\it generalized
simplex\/}.

A {\it Markov kernel\/} between two measurable spaces $\X$ and $\Y$ is a
function $K(x,B)$, where the first argument $x$ ranges over $\X$ and the second
argument ranges over $\B(\Y)$, and such that the following two conditions hold:

\medskip

$\bullet$ $K(\,\cdot\,, B)$ is a measurable function on $\X$  for any fixed
$B\in\B(\Y)$;

$\bullet$ $K(x,\,\cdot\,)$ is a probability measure on $\Y$ for any fixed
$x\in\X$.

\medskip

If both $\X$ and $\Y$ are finite or countably infinite sets, then $K$ is simply
a stochastic matrix of format $\X\times\Y$. About Markov kernels, see, e.g.,
Meyer \cite{Mey66}.

We regard $K$ as a ``link'' between $\X$ and $\Y$ and write this symbolically
as $K:\X\dasharrow\Y$. Sometimes we use the word``link'' as a synonym of
``Markov kernel''. A link is not an ordinary map; this is why we represent it
by a dash arrow. However, it determines a true map $\M(\X)\to\M(\Y)$ taking a
measure $M\in\M(\X)$ to the measure $MK\in\M(\Y)$ defined by
$$
(MK)(B):=\int_{x\in\X}M(dx)K(x,B), \quad B\in\B(\Y).
$$

We regard such a map $\M(\X)\to\M(\Y)$ is an ``affine map'' between
(generalized) simplices. It is a conventional affine map if both $\X$ and $\Y$
are discrete.

Given two links, $K:\X\dasharrow\Y$ and $L:\Y\dasharrow\mathfrak Z$, their
composition $KL:\X\dasharrow\mathfrak Z$ is defined by
$$
(KL)(x,dz)=\int_{y\in\Y}K(x,dy)L(y,dz).
$$
This is a natural generalization of the matrix multiplication.

Thus, one may consider the category whose objects are measurable spaces and
morphisms are links. We need the corresponding notion of {\it projective
limit\/}. To avoid excessive formalism, I define this notion precisely in the
degree of generality that we really need.

Assume we are given an infinite chain of finite or countably infinite sets
together with links between them:
\begin{equation}\label{eq5.C}
\cdots \dasharrow \X_N\dasharrow
\X_{N-1}\dasharrow\cdots\dasharrow\X_2\dasharrow\X_1.
\end{equation}
Because the spaces are discrete, the links are simply stochastic matrices. The
link  between $\X_N$ and $\X_{N-1}$ will be denoted by $\La^N_{N-1}$. The
``categorical'' projective limit of \eqref{eq5.C} is explicitly constructed as
follows.

Chain \eqref{eq5.C} gives rise to a chain of affine maps of simplices
\begin{equation}\label{eq5.D}
\cdots \to \M(\X_N)\to \M(\X_{N-1})\to\cdots\to\M(\X_2)\to\M(\X_1).
\end{equation}
Let $\varprojlim\M(\X_N)$ be the (conventional) projective limit of
\eqref{eq5.D}. By the very definition, elements of $\varprojlim\M(\X_N)$  are
sequences $\{M_N\in\M(\X_N)\}$ such that $M_N\La^N_{N-1}=M_{N-1}$ for every
$N\ge2$. Such a sequence is called a {\it coherent family\/} of measures; here
``coherence'' means that $M_N$'s are consistent with the links. The next
theorem says that $\varprojlim\M(\X_N)$ is a (possibly, generalized) simplex.
More precisely, the claim is the following.

\begin{theorem}\label{thm5.A}
There exists a measurable space $\X_\infty$ and links
$\La^\infty_N:\X_\infty\dasharrow\X_N$, where $N=1,2,\dots$, such that
\begin{equation}\label{eq5.B}
\La^\infty_N\La^N_{N-1}=\La^\infty_{N-1}, \qquad N\ge2,
\end{equation}
and the correspondence $M\mapsto\{M_N:N=1,2,\dots\}$ defined by
$M_N:=M\La^\infty_N$ is a bijection between $\M(\X_\infty)$ and\/
$\varprojlim\M(\X_N)$.

Such a space together with the links $\La^\infty_N$ is unique within a natural
equivalence.
\end{theorem}

A proof based on Choquet's theorem is given in Olshanski \cite[\S9]{Ols03}, a
more general result is contained in Winkler \cite[Chapter 4]{Winkler}.

Note that $\X_\infty$ can be identified with the space of extreme points of the
set $\varprojlim\M(\X_N)$ (which is obviously a convex set), and the nontrivial
part of the theorem is that $\varprojlim\M(\X_N)$ coincides with
$\M(\X_\infty)$.

Note also that for infinite sets $\X_N$ it may happen that
$\varprojlim\M(\X_N)$ is empty (that is, there is no coherent families of
probability measures) and then $\X_\infty$ is empty, too. Here is a simple
example: $\X_N=\{N,N+1,\dots\}$ and the link $\X_N\dasharrow\X_{N-1}$ is
induced by the inclusion $\X_N\subset\X_{N-1}$. However, if all $\X_N$ are
finite sets, then $\X_\infty$ is always nonempty.

Let us regard \eqref{eq5.C} as a kind of discrete time Markov chain with the
transition probabilities determined by the links $\La^N_{N-1}$. (It does not
matter that this chain looks a bit unusual, as the time parameter ranges from
$-\infty$ to 1 and the state space varies with time.) The space
$\varprojlim\M(\X_N)=\M(\X_\infty)$ can be identified with the space of {\it
entrance laws\/} (see Dynkin \cite{Dy-AnnProb78}) for this Markov chain; for
this reason we call $\X_\infty$ the {\it boundary\/} of \eqref{eq5.C}, having
in mind the entrance boundary.

If $\{M_N\}$ is a coherent family of probability measures, then the
corresponding measure $M\in\M(\X_\infty)$ is called the {\it boundary
measure\/} of the family.

By a {\it Markov semigroup\/} on a measurable space $\X$ we mean a semigroup
$P(t)$ of Markov kernels $\X\dasharrow\X$ depending on parameter $t\ge0$. That
is, the kernels are subject to the Chapman--Kolmogorov equation
$P(t_1)P(t_2)=P(t_1+t_2)$ and $P(0)$ is the trivial kernel corresponding to the
identity map $\X\to\X$, that is, $P(0; x,\,\cdot\,)$ is the delta-measure at
$x$.

Under additional assumptions on a Markov semigroup $P(t)$, one can prove that
it serves as the transition function of a Markov process $X$; for instance,
this is so if $P(t)$ is Feller (see Section \ref{sect10} below).

A {\it stationary distribution\/} for a Markov semigroup $P(t)$ is a
probability measure $M\in\M(\X)$ such that  $MP(t)=M$ for all $t\ge0$.

Assume $P(t)$ and $P'(t)$ are Markov semigroups with state spaces $\X$ and
$\X'$, respectively, and $\La:\X\dasharrow\X'$ is a link. We say that $\La$
{\it intertwines\/} the processes if
\begin{equation*}
P(t)\La=\La P'(t), \qquad t\ge0.
\end{equation*}

Now we are in a position to describe a general formalism that we call the {\it
method of intertwiners\/}.

Let us return to the chain \eqref{eq5.C} of discrete spaces and the links
$\La^N_{N-1}:\X_N\dasharrow\X_{N-1}$. Assume that for every $N=1,2,\dots$ we
are given a Markov semigroup $P_N(t)$ on $\X_N$ (that is, simply a semigroup of
stochastic matrices of format $\X_N\times\X_N$), and the links serve as
intertwiners for these semigroups, so that
\begin{equation}\label{eq5.E}
P_N(t)\La^N_{N-1}=\La^N_{N-1}P_{N-1}(t)
\end{equation}
for every $N\ge2$ and any $t\ge0$. We call \eqref{eq5.E} the {\it master
relation\/}. Finally, assume that the boundary $\X_\infty$ of \eqref{eq5.C} is
nonempty.

\begin{theorem}\label{thm5.B}
{\rm(i)} Under these hypotheses there exists a unique Markov semigroup
$P_\infty(t)$ on $\X_\infty$ such that
$$
P_\infty(t)\La^\infty_N=\La^\infty_N P_N(t), \qquad N=1,2,\dots, \quad t\ge0.
$$

{\rm(ii)} Assume additionally that there exists a coherent family $\{M_N\}$ of
probability distribution such that $M_N$ is is a stationary distribution for
$P_N(t)$ for every $N$. Then the corresponding boundary measure on $\X_\infty$
is a stationary distribution for $P_\infty(t)$.
\end{theorem}

These assertions are direct consequences of the definitions. We call
$P_\infty(t)$ the {\it boundary Markov semigroup\/}.

It may well happen that the semigroups $P_N(t)$ are not given in an  explicit
form. Then, to check the master relation \eqref{eq5.E}, one may try to reduce
it to its infinitesimal version,
\begin{equation}\label{eq5.F}
D_N\La^N_{N-1}=\La^N_{N-1}D_{N-1},
\end{equation}
where $D_N$ stands for the infinitesimal generator of $P_N(t)$ and
\eqref{eq5.F} should be understood as a relation for operators acting in
suitable function spaces (see Section \ref{sect10} below); when applied to a
function $F$, \eqref{eq5.F} should be read from right to left:
\begin{equation*}
D_N\La^N_{N-1}F=\La^N_{N-1}D_{N-1}F.
\end{equation*}

\section{Examples}\label{sect6}

Here I illustrate the formalism of the preceding section by two simple
examples.

Consider the following chain of type \eqref{eq5.C} coming from the {\it Pascal
triangle\/}: the spaces are finite sets,
$$
\X_N=\{0,1,\dots,N\}\subset\Z,
$$
and the links are defined by
\begin{equation*}
\La^N_{N-1}(n,m)=\begin{cases} \dfrac{N-n}N, & m=n,\\ \dfrac nN, & m=n-1,\\ 0,
& m\ne n,n-1. \end{cases}
\end{equation*}

The boundary $\X_\infty$ of this chain can be identified with the closed unit
interval $[0,1]$ and the links $\La^\infty_N$ are given by
\begin{equation*}
\La^\infty_N(x,n)=\binom{N}{n}x^n(1-x)^{N-n}, \qquad x\in[0,1], \quad
n=0,1,\dots,N;
\end{equation*}
that is, $\La^\infty_N(x,\,\cdot\,)$ is the binomial distribution with
parameter $x$. This fact is equivalent to de Finetti theorem or else to the
solution of the Hausdorff moment problem. See, e.g., Gnedin--Pitman
\cite{GnedinPitman} and references therein.

\begin{example}\label{ex6.A}

Fix two real parameters $a>0$, $b>0$. We are going to define, for every $N$, a
continuous time Markov chain $X_N$ on $\X_N$. To do this we exhibit its
generator $D_N$, which is a difference operator on $\X_N\subset\Z$. Its action
on a test function $F$ is given by
\begin{equation}\label{eq6.A}
\begin{aligned}
(D_NF)(n)=(N-n)(n+a)&\left[F(n+1)-F(n)\right]\\+n(N+b-n)&\left[F(n-1)-F(n)\right].
\end{aligned}
\end{equation}

It is directly verified that the operators $D_N$ satisfy \eqref{eq5.F}, from
which one can deduce that the corresponding semigroups $P_N(t)$ (the transition
functions of the chains $X_N$) satisfy the master relation \eqref{eq5.E}.

Therefore, by virtue of Theorem \ref{thm5.B}, part (i), these semigroups give
rise to a boundary Markov semigroup $P_\infty(t)$ on $[0,1]$. One can prove
that $P_\infty(t)$ is the transition function of a diffusion process $X_\infty$
on $[0,1]$ with the infinitesimal  generator
\begin{equation*}
D_\infty=x(1-x)\frac{d^2}{dx^2}+[a-(a+b)x]\frac{d}{dx}.
\end{equation*}

For every $N$, the chain $X_N$ has a unique stationary distribution $M_N$,
$$
M_N(n)=\frac{\Ga(a+b)N!}{\Ga(a)\Ga(b)\Ga(a+b+N)}\,
\frac{\Ga(a+n)\Ga(b+N-n)}{n!(N-n)!}\,, \qquad n=0,1,\dots,N,
$$
which is the well-known hypergeometric distribution (see, e.g., Feller
\cite[ch. II, \S6]{Feller-Book}). The sequence $\{M_N\}$ is a coherent family,
and the corresponding boundary measure $M_\infty$ is Euler's beta distribution,
$$
M_\infty(dx)=\frac{\Ga(a+b)}{\Ga(a)\Ga(b)}x^{a-1}(1-x)^{b-1}dx, \qquad
x\in[0,1].
$$
It is a unique stationary distribution for $X_\infty$.
\end{example}

\begin{example}\label{ex6.B}
Here the sets $\X_N$ and the links $\La^N_{N-1}$ are as in the preceding
example, but we choose different Markov chains: this time they are defined by
the difference operators
\begin{equation}\label{eq6.B}
\begin{aligned}
(D'_NF)(n)=(N-n)c&\left[F(n+1)-F(n)\right]\\+n(1-c)&\left[F(n-1)-F(n)\right],
\end{aligned}
\end{equation}
where $c\in[0,1]$ is a fixed parameter. The corresponding Markov semigroups
$P'_N(t)$ are again consistent with the links and so give rise to a boundary
Markov semigroup $P'_\infty(t)$ on $[0,1]$. But $P'_\infty(t)$ turns out to be
degenerate in the sense that it corresponds to a deterministic process
$X'_\infty$. Namely, the generator of $X'_\infty$ is a first order differential
operator,
$$
D'_\infty=(c-x)\frac{d}{dx},
$$
so that $X'_\infty$ is not a genuine Markov process but a flow of endomorphisms
of the interval $[0,1]$ generated by a vector field. The dynamics is easily
described: the (deterministic) trajectory $x(t)$ issued from a given point
$x(0)$ has the form
$$
x(t)=c-(c-x(0))e^{-t}, \qquad t\ge0,
$$
so that the interval $[0,1]$ is contracted to the point $c$ exponentially fast.

Note that the $N$th chain $X'_N$ has a unique stationary distribution $M'_N$,
which is the binomial distribution with parameter $c$,
$$
M'_N(n)=\binom Nn c^n(1-c)^{N-n}, \qquad n=0,1,\dots, N.
$$
The distributions $M'_N$ form a coherent system with the boundary measure
$\de_c$, the delta-measure at $c\in[0,1]$. This agrees with the evident fact
that $\de_c$ is a (unique) stationary distribution of the flow $X'_\infty$.
\end{example}

As seen from the second example, it may happen that a boundary process
constructed according to the general scheme of Section \ref{sect5} degenerates
to a deterministic process. So, if one wants to get a genuine Markov dynamics
on the boundary, one needs additional arguments guaranteeing that such a
degeneration does not occur.

\section{Extremal characters of $U(\infty)$ and the boundary $\Om_\infty$}
\label{sect7}

Here I introduce certain links $\La^N_{N-1}:\Om_N\dasharrow\Om_{N-1}$ between
the subsets $\Om_N$ (see their definition in Section \ref{sect4}) and discuss
the meaning of the boundary $\Om_\infty$ of the chain
\begin{equation}\label{eq7.A}
\cdots \dasharrow \Om_N\dasharrow
\Om_{N-1}\dasharrow\cdots\dasharrow\Om_2\dasharrow\Om_1.
\end{equation}
In the end of Section \ref{sect4}, it was pointed out that $\Om_N$
parameterizes the dual object $\wh{U(N)}$. I will explain this point in more
detail.

By definition, the {\it dual object\/} $\wh G$ to a compact group $G$ is the
set of equivalence classes of irreducible finite-dimensional representations of
$G$. As well known, the irreducible representations of the group $G=U(N)$ are
indexed by the vectors $\la=(\la_1,\dots,\la_N)\in\Z^N$ with nonincreasing
coordinates, $\la_1\ge\dots\ge\la_N$; such vectors are called {\it signatures
of length\/} $N$ (see, e.g., Weyl \cite{Wey39}, Zhelobenko \cite{Zhe70}). There
is a one-to-one correspondence $\la\leftrightarrow\ell$ between signatures
$\la$ and elements $\ell\in\Om_N$ given by
$$
l_i=\la_i+N-i, \qquad i=1,\dots,N.
$$
Thus, we may take $\Om_N$ as the set of parameters for the dual object to
$U(N)$.

Besides the parameterization of $\wh{U(N)}$, the only extra fact about
representations that we need is the Gelfand--Tsetlin {\it branching rule\/}
which describes the decomposition of an irreducible representation of $U(N)$
when restricted to $U(N-1)\subset U(N)$. Here $U(N-1)$ is considered as the
subgroup of $U(N)$ that fixes the last basis vector in $\C^N$.

Let us introduce some notation. The irreducible representation of $U(N)$
corresponding to a signature $\la$ will be denoted by $\pi^{\la,N}$. Two
signatures $\la=(\la_1,\dots,\la_N)$ and $\mu=(\mu_1,\dots,\mu_{N-1})$ are said
to be {\it interlaced\/} if
$$
\la_1\ge\mu_1\ge\la_2\ge\dots\ge\la_{N-1}\ge\mu_{N-1}\ge\la_N;
$$
then we write $\mu\prec\la$.

The branching rule (Gelfand--Tsetlin \cite{GT}, Zhelobenko \cite{Zhe70}) says
that
\begin{equation*}
\pi^{\la,N}\big|_{U(n-1)}=\bigoplus_{\mu:\, \mu\prec\la}\pi^{\mu,N-1}.
\end{equation*}
Taking the dimensions of the both sides gives the identity
\begin{equation*}
\dim\pi^{\la,N}\big|_{U(n-1)}=\sum_{\mu:\, \mu\prec\la}\dim\pi^{\mu,N-1}.
\end{equation*}
We use it to define a link $\La^N_{N-1}:\Om_N\dasharrow\Om_{N-1}$, as follows.
Let $\ell\in\Om_N$, $\ell'\in\Om_{N-1}$, and let $\la\leftrightarrow\ell$ and
$\mu\leftrightarrow\ell'$ be the corresponding signatures. We set
\begin{equation*}
\La^N_{N-1}(\ell,\ell'):=\begin{cases}\dfrac{\dim\pi^{\mu,N-1}}{\dim\pi^{\la,N}},
& \mu\prec\la,\\ 0, & \textrm{otherwise.} \end{cases}
\end{equation*}
Because of the above identity, $\La^N_{N-1}$ is a stochastic matrix, so the
definition is correct. Thus, we have constructed the chain \eqref{eq7.A}.

Observe that the branching rule entails a direct combinatorial definition of
the quantity $\dim\pi^{\la,N}$: namely, it is equal to the total number of
sequences
\begin{equation*}
\la^{(1)}\prec\la^{(2)}\prec\dots\prec\la^{(N)}=\la,
\end{equation*}
where $\la^{(i)}$ is a signature of length $i$. Such sequences are often
written as triangular arrays, called {\it Gelfand--Tsetlin schemes\/} or {\it
Gelfand--Tsetlin patterns\/}, see Gelfand--Tsetlin \cite{GT}, Zhelobenko
\cite{Zhe70}. On the other hand, there is an explicit formula, which is a
particular case of {\it Weyl's dimension formula\/} (Weyl \cite{Wey39},
Zhelobenko \cite{Zhe70}):
\begin{equation*}
\dim\pi^{\la,N}=\prod_{1\le i<j\le
N}\frac{\la_i-\la_j+j-i}{j-i}=\frac{\prod_{1\le i<j\le
N}(\ell_i-\ell_j)}{1!\dots(N-1)!}.
\end{equation*}
This makes the definition of the links formally independent of the
representation theory of the unitary groups.

Let $\Om_\infty$ stand for the boundary of the chain \eqref{eq7.A}. Observe
that the chain of Section \ref{sect6} can be embedded into \eqref{eq7.A}.
Namely, the element $n\in\X_N$ is identified with the signature of length $N$
of the form $(1,\dots,1,0,\dots,0)$, where the number of 1's equals $n$. This
shows that the boundary $\Om_\infty$ contains the boundary of the chain of
Section \ref{sect6}. In particular, it follows that $\Om_\infty$ is nonempty.

Let us define the group $U(\infty)$ as the union of the groups $U(N)$ embedded
one into another as indicated above; $U(\infty)$ belongs to the class of {\it
inductive limits of compact groups\/}. We extend the definition of dual object
to groups $G=\varinjlim G_N$ from this class in the following way.

A function $\chi:G\to\C$ is said to be an {\it extremal character\/} if it
satisfies the following three conditions:

\medskip
$\bullet$ First, $\chi$ is {\it normalized\/}, that is, $\chi(e)=1$.

\smallskip

$\bullet$ Second, $\chi$ is {\it central\/} meaning that it is constant on each
conjugacy class.

\smallskip

$\bullet$ Third, for any elements $g_1, g_2\in G$ one has
\begin{equation}\label{eq7.B}
\lim_{N\to\infty}\int_{h\in G_N}
\chi(g_1hg_2h^{-1})m_{G_N}(dh)=\chi(g_1)\chi(g_2),
\end{equation}
where $m_{G_N}$ denotes the normalized Haar measure on the compact group $G_N$.
We define $\wh G$ as the set of all such functions.

\medskip

Here is an explanation why this new definition of $\wh G$ extends the previous
one. For a compact group $G$, every irreducible representation is uniquely
determined by its {\it character\/} $\chi^\pi$: this is a function on $G$ given
by
\begin{equation*}
\chi^\pi(g)=\tr(\pi(g)), \qquad g\in G.
\end{equation*}
The normalized function
$$
\wt\chi^{\,\pi}(g):=\chi^\pi(g)/\chi^\pi(e)=\chi^\pi(g)/\dim \chi^\pi
$$
is called a {\it normalized irreducible character\/}. By the very definition,
such functions may serve as parameters for the dual object $\wh G$. On the
other hand, it is well known that the normalized irreducible characters of a
compact group $G$ are precisely those functions $\chi: G\to\C$ that satisfy the
first and second conditions stated above and the simplified form of the third
condition, the {\it functional equation\/}
\begin{equation}\label{eq7.C}
\int_{h\in G} \chi(g_1hg_2h^{-1})m_G(dh)=\chi(g_1)\chi(g_2), \qquad \forall
g_1,g_2\in G
\end{equation}
(here $m_G$ is the normalized Haar measure on $G$).

There is another but equivalent definition of the extremal characters which
shows that they are extreme points of a certain convex set of functions on $G$,
see Olshanski \cite[\S1]{Ols-AA89}, \cite[\S\S23-24]{Ols-G&B}. These references
also explain how extremal characters are related to unitary representations.

\begin{theorem}\label{thm7.A}
There exists a natural one-to-one bijective correspondence
$\Om_\infty\leftrightarrow\wh{U(\infty)}$.
\end{theorem}

Thus, the boundary $\Om_\infty$ has a representation-theoretic meaning. Note,
however, that this theorem is a purely abstract result that provides no
information about the size of the boundary. Its explicit description is given
below in Section \ref{sect9}.

\section{The master equation and the stationary distribution}\label{sect8}

Recall that in Section \ref{sect4} we constructed Markov chains $\wt X^\uv_N$
on the sets $\Om_N$. In Section \ref{sect4}, parameters $N$ and $\uv$ were
fixed, while in what follows $N$ will vary and $\uv$ will vary together with
$N$. Namely, let us set
\begin{equation}\label{eq8.B}
\uv=:=(z+N-1,z'+N-1, w,w'),
\end{equation}
where $\zw$ is a fixed quadruple of parameters such that $(z,z')$ and $(w,w')$
are admissible. According to this we slightly change the notation. Let us set
\begin{equation*}
X^\zw_N=\wt X^{(z+N-1,z'+N-1,w,w')}
\end{equation*}
and denote by $P^\zw_N(t)$ the Markov semigroup of $X^\zw_N$.

\begin{theorem}\label{thm8.A}
The semigroups just defined and the links
$\La^N_{N-1}:\Om_N\dasharrow\Om_{N-1}$ defined in Section \ref{sect7} satisfy
the master equation \eqref{eq5.E}. That is,
\begin{equation}\label{eq8.A}
P^\zw_N(t)\La^N_{N-1}=\La^N_{N-1}P^\zw_{N-1}(t).
\end{equation}
\end{theorem}

By virtue of Theorem \ref{thm5.B}, part (i), the semigroups $P^\zw_N(t)$ give
rise to a boundary semigroup on $\Om_\infty$; let us denote it by
$P^\zw_\infty(t)$.

Set $M^\zw_N:=\wt M^\uv_N$, where $\wt M^\uv_M$ is the invariant measure from
Theorem \ref{thm4.A}, part (ii), defined by \eqref{eq4.D}; here, as above, the
quadruple $\uv$ is given by \eqref{eq8.B}. Let us assume additionally that
$z+z'+w+w'>-1$. Then the constant factor in \eqref{eq4.D} can be chosen so that
$M^\zw_N$ becomes a probability measure.

\begin{theorem}\label{thm8.B}
Assume $z+z'+w+w'>-1$. The probability measures $M^\zw_N$ just defined satisfy
the relation
\begin{equation}\label{eq8.C}
M^\zw_N\La^N_{N-1}=M^\zw_{N-1}, \qquad N\ge2,
\end{equation}
so that $\{M^\zw_N: N\ge1\}$ is a coherent family that determines a boundary
measure $M^\zw$ on\/ $\Om_\infty$.

This measure is a stationary distribution for the boundary semigroup
$P^\zw_\infty(t)$.
\end{theorem}

The first assertion of the theorem was proved in Olshanski \cite{Ols03}. The
second assertion is a formal consequence of the first assertion and the fact
that $M^\zw_N$ is a stationary distribution for $P_N(t)$ for every $N$.

The measures $M^\zw_\infty$ are called the {\it boundary zw-measures\/}. Note
that $M^\zw_\infty$ does not change under transposition $z\leftrightarrow z'$
or $w\leftrightarrow w'$.

\begin{theorem}[Gorin \cite{Gorin-FAA}]\label{thm8.C}
The zw-measures corresponding to different, up to the above transpositions,
quadruples of parameters are pairwise disjoint, that is, mutually singular.
\end{theorem}

\section{The Edrei--Voiculescu theorem}\label{sect9}

Consider a two-sided infinite sequence $\{\varphi_n: n\in\Z\}$ of real numbers
and assign to it the two-sided infinite Toeplitz matrix $T$ with the entries
$T(i,j):=\varphi_{j-i}$, where $i,j\in\Z$. The sequence $\{\varphi_n\}$ is
called {\it totally positive\/} if all minors of $T$ are nonnegative. In more
detail, the minors of order 1 are the numbers $\varphi_n$, so they  must be
nonnegative; next, the minors of order 2 are indexed by two arbitrary couples
of integers, $n_1<n_2$ and $m_1<m_2$, which leads to the condition
$\varphi_{n_1}\varphi_{m_2}-\varphi_{n_2}\varphi_{m_1}\ge0$, and so on.

The problem of classification of the totally positive sequences was posed by
Schoenberg and solved by Edrei \cite{Edrei}. The result is deep and the answer
is beautiful. To avoid excessive and unnecessary complication we will impose
the additional requirement that
\begin{equation}\label{eq9.B}
\sum_{n\in\Z}\varphi_n=1.
\end{equation}
Then the result is conveniently stated in terms of the generating series
\begin{equation}\label{eq9.A}
\Phi(u):=\sum_{n\in\Z}\varphi_n u^n.
\end{equation}
Because of \eqref{eq9.B}, the series converges on the unit circle $\T\subset\C$
and represents there a continuous function.

\begin{theorem}[Edrei \cite{Edrei}]\label{thm9.A}
The totally positive sequences with the normalization condition \eqref{eq9.B}
are parameterized by sextuples $\om=(\al^+,\be^+,\al^-,\be^-,\de^+,\de^-)$,
where $\al^\pm$ and $\be^\pm$ are infinite sequences of nonincreasing
nonnegative reals $\{\al^\pm_i: i=1,2,\dots\}$ and $\{\be^\pm_i:
i=1,2,\dots\}$, respectively, such that
$$
\sum_{i=1}^\infty (\al^\pm_i+\be^\pm_i)\le\de^\pm, \qquad \be^+_1+\be^-_1\le1.
$$

Given such a sextuple $\om$, the generating series of the corresponding
sequence has the form
\begin{equation}\label{eq9.C}
\Phi(u;\om)=e^{\ga^+(u-1)+\ga^-(u^{-1}-1)} \prod_{i=1}^\infty
\frac{(1+\be^+_i(u-1))}{(1-\al^+(u-1))}\frac{(1+\be^-_i(u^{-1}-1))}{(1-\al^-(u^{-1}-1))}\,,
\end{equation}
where
$$
\ga^\pm:=\de^\pm-\sum_{i=1}^\infty (\al^\pm_i+\be^\pm_i).
$$
\end{theorem}

In particular, the generating series converges in an annulus around $\T$ and
extends to a meromorphic function in $\C$. (About the theory of total
positivity see Karlin's fundamental monograph \cite{Karlin}.)

Voiculescu discovered that the same functions \eqref{eq9.C} appear in the
context of the representation theory of the group $U(\infty)$. Namely, the
following result holds. (Below we use the fact that every matrix $U\in
U(\infty)$ is conjugated to a diagonal matrix whose diagonal entries, the
eigenvalues of $U$, lie on the unit circle $\T$ and only finitely many of them
are distinct from 1.)

\begin{theorem}\label{thm9.B}
The extremal characters of $U(\infty)$ are precisely the functions of the form
\begin{equation}\label{eq9.D}
\chi_\om(U)=\prod_{j=1}^\infty \Phi(u_j;\om), \quad U\in U(\infty),
\end{equation}
where $\om$ ranges over the same collection of parameters as in Theorem
\ref{thm9.A}, and $u_1,u_2,\dots$ are the eigenvalues of the matrix $U$.
{\rm(}Note that the infinite product here is actually finite, because
$\Phi(1;\om)=1$ and $u_j=1$ for $j$ large enough.{\rm)}
\end{theorem}

Thus, the extremal characters of $U(\infty)$ and the totally positive sequences
are in one-to-one correspondence, so that the classification problems for these
two kinds of objects coincide.

\begin{remark}
Here are brief historical comments concerning  Theorem \ref{thm9.B}. Voiculescu
was the first person to study the extremal characters of $U(\infty)$ (see his
paper \cite{Voi76}). He proved that all functions of the form \eqref{eq9.D} are
extremal characters. He also explained why the extremal characters should be
given by multiplicative expressions with respect to the eigenvalues. He did not
prove that the list of Theorem \ref{thm9.B} is exhaustive, but obtained some
partial results in this direction. Then Vershik--Kerov \cite{VK82} and Boyer
\cite{Boyer} independently drew attention to the earlier work of Edrei
\cite{Edrei}, of which Voiculescu was unaware. Boyer explained how to deduce
Theorem \ref{thm9.B} from Edrei's theorem. Vershik and Kerov sketched quite a
different approach to Theorem \ref{thm9.B}, already tested on the example of
the infinite symmetric group \cite{VK81}. A detailed proof (in a broader
context), based on the ideas of \cite{VK82}, appeared later in
Okounkov--Olshanski \cite{OO98}. Recently, one more proof was proposed in
Borodin--Olshanski \cite{BO-GT-Appr}, and soon after that Petrov \cite{Petrov}
found a simpler version of it together with a generalization.
\end{remark}

\begin{remark}
It is worth noting that the multiplicativity property of extremal characters of
$U(\infty)$ is related to specific properties of some infinite-dimensional
groups and does not hold for finite-dimensional (noncommutative) groups. The
nature of this phenomenon is analyzed in my expository paper
\cite{Ols-Semigroups} (see also \cite{Ols-G&B}). One of the explanations given
in \cite{Ols-Semigroups} is related to a concentration property for the Haar
measure of $U(N)$ (and other similar groups) as $N\to\infty$. Although the
normalized irreducible characters of the groups $U(N)$ are not multiplicative,
they become ``approximately multiplicative'' as $N$ gets large. This can be
seen from \cite{Ols-Semigroups}, and recently, Gorin and Panova \cite{GP13}
found new character formulas which demonstrate this effect in a very clear
manner.
\end{remark}

Theorem \ref{thm9.B} shows that the boundary $\Om_\infty$, whose abstract
definition was given in section \ref{sect8}, admits an explicit description.
Namely, it can be identified with the region in the product space
$$
\R^{4\infty+2}:=\R^\infty\times\R^\infty\times\R^\infty\times\R^\infty\times\R\times\R
$$
formed by the sextuples $\om=(\al^+,\be^+,\al^-,\be^-,\de^+,\de^-)$ from
Theorem \ref{thm9.A}.

To complete the description of the boundary it remains to specify the links
$\La^\infty_N:\Om_\infty\dasharrow \Om_N$:
\begin{equation}\label{eq9.E}
\La^\infty_N(\om,\ell)=\dim\pi^{\la,N}\det\left[\varphi_{\la_i-i+j}\right]_{i,j=1}^N,
\end{equation}
where $\om\in\Om_\infty\subset\R^{4\infty+2}$, $\ell$ ranges over $\Om_N$,
$\la$ is the signature corresponding to $\ell$, and $\{\varphi_n\}$ is the
collection of the Laurent coefficients of the function $\Phi_\om(u)$ defined by
\eqref{eq9.C}. Note that the right-hand of \eqref{eq9.E} is nonnegative (as it
should be), because the determinant $\det\left[\varphi_{\la_i-i+j}\right]$ is
nonnegative due to the total positivity of $\{\varphi_n\}$.

\begin{remark}
The determinants appearing in \eqref{eq9.E} do not exhaust all minors of the
Toeplitz matrix $T$ associated with the sequence $\{\varphi_n\}$: indeed, these
are minors with consecutive column numbers. However, the nonnegativity of these
special minors already suffices to conclude that all minors of $T$ are
nonnegative, see Boyer \cite{Boyer}. This fine point is necessary for
establishing the bijection between totally positive sequences and extremal
characters.
\end{remark}

\begin{remark}
Now one can explain how the boundary $[0,1]$ of the Pascal triangle discussed
in Section \ref{sect7} is located in $\Om_\infty$. Namely, the interval $[0,1]$
is identified with the set of those $\om$'s for which $\be^+_1=\de^+=x\in[0,1]$
and all other coordinates of $\om$ equal 0.
\end{remark}

\section{The generator}\label{sect10}

I will start with a few definitions and facts concerning Feller Markov
processes. For more detail, see Ethier--Kurtz \cite{EK86}, Liggett \cite{Lig}.

Assume $\X$ is a locally compact metrizable separable space and denote by
$C_0(\X)$ the space of real-valued continuous functions vanishing at infinity.
Let us equip $C_0(\X)$ with the supremum norm; then it becomes a separable
Banach space. Note that the larger Banach space $C(\X)$ of bounded continuous
functions is not separable unless $\X$ is compact, in which case
$C_0(\X)=C(\X)$; but we are interested in the case when $\X$ is not compact.

A semigroup $P(t)$ of Markov kernels on $\X$ is said to be a {\it Feller
semigroup\/} if it preserves the space $C_0(\X)$ and induces in it a strongly
continuous operator semigroup. Then $P(t)$ generates a Markov process $X$ on
$\X$ with sufficiently good sample trajectories; $X$ is called a {\it Feller
process\/}.

A Feller semigroup is uniquely determined by its generator $A$, which is a
closed dissipative operator on $C_0(\X)$. Its domain $\Dom(A)$ is  formed by
those elements $F\in C_0(\X)$ for which the limit
$$
AF:=\lim_{t\to0}\frac{P(t)F-F}t
$$
exists. A  subspace $\mathscr F\subset\Dom(A)$ is called a {\it core\/} for $A$
if the closure of the restriction of $A$ to $\mathscr F$ coincides with $A$. In
practice, it is usually problematic to explicitly describe $\Dom(A)$, and then
one is satisfied by indicating the action of $A$ on an appropriate core,
because this suffices to specify $A$.

Now let us return to our boundary semigroup $P^\zw_\infty(t)$. An important
fact is that the boundary $\Om_\infty$ is a locally compact space with respect
to the topology induced by the product topology of the ambient infinite product
space $\R^{4\infty+2}$.

\begin{theorem}\label{thm10.A}
$P^\zw_\infty(t)$ is a Feller semigroup, so it give rise to a Feller process
$X^\zw_\infty$ on $\Om_\infty$.
\end{theorem}

This result raises the question about the semigroup generator as an operator in
the Banach space $C_0(\Om_\infty)$. One can exhibit a core $\mathscr F\subset
C_0(\Om_\infty)$ and prove that the action of the generator on $\mathscr F$ is
implemented by a second order differential operator $D^\zw_\infty$ with
countably many variables.

At first glance, one would expect that $D^\zw_\infty$ is somehow written in
terms of the natural coordinate system $(\al^\pm_i; \be^\pm_i; \de^\pm)$ on
$\Om_\infty$, but it is not so. The natural coordinates are unsuitable, and we
have to pass to other variables that are (in some sense) supersymmetric
functions of the natural coordinates. These new variables are the Laurent
coefficients of \eqref{eq9.C}, which we denoted by $\varphi_n$, $n\in\Z$. As
the core $\mathscr F$ we take a certain subspace in $\R[\ldots,\varphi_{-1},
\varphi_0,\varphi_1,\ldots]$, the algebra of polynomials in variables
$\varphi_n$. Then $D^\zw_\infty$ is written in the form
\begin{equation}\label{eq10.A}
D^\zw_\infty=\sum_{n_1,n_2\in\Z}\Ga^{(2)}_{n_1,n_2}\frac{\pd^2}{\pd
\varphi_{n_1}\varphi_{n_2}}+\sum_{m\in\Z}\Ga^{(1)}_m \frac{\pd}{\pd \varphi_m},
\end{equation}
where the coefficients $\Ga^{(2)}_{n_1,n_2}$ are certain infinite quadratic
expressions in variables $\varphi_n$ while the coefficients $\Ga^{(1)}_m$ are
certain finite linear combinations of these variables.

Note that only coefficients $\Ga^{(1)}_m$ depend on the basic parameters $\zw$
while coefficients $\Ga^{(2)}_{n_1,n_2}$ do not. This implies the following.
Recall that the boundary process admits a stationary distribution, the
zw-measure $M^\zw$ (see Theorem \ref{thm8.B}). Consider the Hilbert space
$H^\zw:=L^2(\Om_\infty,M^\zw)$ and introduce the (pre)Dirichlet form
corresponding to $D^\zw_\infty$,
$$
\mathscr E(F,G):=-(D^\zw_\infty F,G), \qquad F,G\in\mathscr F;
$$
here the brackets in the right-hand side denote the inner product in $H^\zw$.
Then we get
$$
\mathscr E(F,G)=\int_{\om\in\Om_\infty} \Ga(F,G)M^\zw(d\om), \quad
\Ga(F,G):=\sum_{n_1,n_2\in\Z}\Ga^{(2)}_{n_1,n_2} \frac{\pd F}{\pd
\varphi_{n_1}}\frac{\pd G}{\pd \varphi_{n_2}},
$$
where only $M^\zw$ depends on the basic parameters while the form $\Ga(F,G)$
does not. One may speculate that this form somehow expresses the ``inner
geometry'' of the space $\Om_\infty$.

The fact that $D^\zw_\infty$ has second order implies that $X^\zw_\infty$
cannot degenerate to a deterministic process, as in Example \ref{ex6.B}. This
conclusion can be also deduced from the fact that $M^\zw_\infty$ is not only a
stationary distribution but also a symmetrizing measure.

\section{Summary}\label{sect11}

The starting point of the story is a 4-parameter family
$\{X^\zw_N:N=1,2,\dots\}$ of continuous time Markov chains on the dual objects
$\Om_N=\wh{U(N)}$. For any fixed quadruple $\zw$ of parameters, the chains
$X^\zw_N$ are consistent with some canonical ``links'' (stochastic matrices)
relating the sets $\Om_N$ to each other. This makes it possible to apply the
abstract ``method of intertwiners'' and establish the existence of Markov
semigroups $P^\zw_\infty(t)$ on an infinite-dimensional locally compact space
$\Om_\infty=\wh{U(\infty)}$. Every semigroup $P^\zw(t)$ possesses the Feller
property and so determines a Feller Markov process $X^\zw_\infty$ on
$\Om_\infty$. This process has a unique stationary distribution $M^\zw_\infty$,
which also serves as a symmetrizing measure. The action of the infinitesimal
generator of the process on an appropriate core can be explicitly described,
and it turns out that it is implemented by a second order differential operator
with infinitely many variables.

\section{Concluding remarks}

The boundary zw-measures $M^\zw_\infty$ are of great interest for harmonic
analysis (Olshanski \cite{Ols03}). However, they are defined indirectly,
through an abstract existence theorem, which makes it difficult to work with
them. As seen from Theorem \ref{thm8.C}, the zw-measures cannot be given by
densities with respect to a reference measure on $\Om_\infty$.

A way to describe the zw-measures is to interpret them as the laws of some {\it
determinantal point processes\/} whose correlation kernels can be explicitly
computed (Borodin--Olshanski \cite{BO-AnnMath}), so every zw-measure is a {\it
determinantal measures\/}. (About such measures, see Borodin
\cite{Bor-Oxford11} and references therein.)

The above results show that $M^\zw_\infty$ can also be characterized as the
only invariant measure of the process $X^\zw_\infty$.

In view of the results of \cite{BO-AnnMath}, it seems plausible that the
multi-time finite-dimensional distributions of the process $X^\zw_\infty$
started from corresponding the zw-measure also have the determinantal
structure.

It would be very interesting to learn more about the properties of the
processes $X^\zw_\infty$. For instance, is it true that the sample trajectories
of $X^\zw_\infty$ are continuous?

The present notes do not cover all the results of the paper Borodin--Olshanski
\cite{BO-GT-Dyn}. As shown in that paper (see also \cite{BO-AnnMath}),
$X^\zw_\infty$ can be interpreted as a time-dependent point process with
infinitely many particles. Although the interaction between the particles is
highly nonlocal, it turns out that $X^\zw_\infty$ can be obtained as a
projection of another Markov process, in which the interaction between the
particles is local (see \cite[\S9]{BO-GT-Dyn}).

Finally, note that there exists a parallel theory in which the role of
$U(\infty)$ is played by the infinite symmetric group, see Borodin--Olshanski
\cite{BO-MMJ}, \cite{BO-ThomaCone} and references therein.


\begin{thebibliography}{Macd95}

\bibitem{Bor-Oxford11}
A. Borodin, {\it Determinantal point processes\/}. In: {\it The Oxford Handbook
on Random Matrix Theory\/}, Gernot Akemann, Jinho Baik, and Philippe Di
Francesco, eds.  Oxford University Press, 2011, Chapter 11, 231--249;
arXiv:0911.1153.

\bibitem{BG-PTRF}
A. Borodin and V. Gorin, {\it Markov processes of infinitely many
nonintersecting random walks\/}. Probab. Theory  Rel.  Fields {\bf155} (2013),
Issue 3-4, 935--997; arXiv:1106.1299.

\bibitem{BG-Lectures}
A. Borodin and V. Gorin, {\it Lectures on integrable probability\/}. Preprint,
arXiv:1212.3351

\bibitem{BO-AnnMath}
A. Borodin and G. Olshanski, {\it Harmonic analysis on the infinite-dimensional
unitary group and determinantal point processes\/}. Ann. Math. {\bf161} (2005),
1--104; arXiv:math/0109194.

\bibitem{BO05-EurCongr}
A. Borodin and G. Olshanski, {\it Representation theory and random point
processes\/}. In: {\it European Congress of Mathematics\/}, Eur. Math. Soc.,
Zurich, 2005, pp. 73--94; arXiv:math/0409333.


\bibitem{BO-PTRF}
A.~Borodin and G.~Olshanski, {\it Markov processes on partitions\/}. Probab.
Theory  Rel. Fields {\bf135} (2006), 84--152; arXiv:math-ph/0409075.


\bibitem{BO-GT-Dyn}
A. Borodin and G. Olshanski, {\it Markov processes on the path space of the
Gelfand-Tsetlin graph and on its boundary\/}  J. Funct. Anal.{\bf263} (2012),
248--303; arXiv:1009.2029.

\bibitem{BO-GT-Appr}
A.~Borodin and G.~Olshanski, {\it The boundary of the Gelfand--Tsetlin graph: A
new approach\/}. Advances in Math. {\bf230} (2012), 1738--1779;
arXiv:1109.1412.

\bibitem{BO-MMJ}
A.~Borodin and G.~Olshanski, {\it The Young bouquet and its boundary\/}. Moscow
Math. J. {\bf13} (2013), no. 2, 191--230; arXiv:1110.4458.

\bibitem{BO-ThomaCone}
A.~Borodin and G.~Olshanski, {\it Markov dynamics on the Thoma cone: a model of
time-dependent determinantal processes with infinitely many particles\/}.
Preprint, arXiv:1303.2794.

\bibitem{Boyer}
R. P. Boyer,  {\it Infinite traces of AF-algebras and characters of
$U(\infty)$\/}. J. Operator Theory {\bf9} (1983), 205--236.

\bibitem{Dougall}
J.~Dougall, {\it On Vandermonde's theorem and some more general expansions\/}.
Proc. Edinburgh Math. Soc. {\bf 25} (1907),  114--132.

\bibitem{Dy-AnnProb78}
E. B. Dynkin, {\it Sufficient statistics and extreme points\/}. Ann. Probab.
{\bf6} (1978), 705--730.

\bibitem{Dys}
F. J. Dyson, {\it A Brownian-motion model for the eigenvalues of a random
matrix\/}. J. Math. Phys. {\bf3} (1962) 1191--1198.

\bibitem{Er}
A.~Erdelyi (ed.) {\it Higher transcendental functions\/}, {\rm I}.
McGraw--Hill, New York, 1953.

\bibitem{Edrei}
A. Edrei, {\it On the generating function of a doubly infinite, totally
positive sequence\/}. Trans. Amer. Math. Soc. {\bf74} (1953), 367--383.

\bibitem{EK86}
S. N. Ethier and T. G. Kurtz, {\it Markov processes --- Characterization and
convergence\/}. Wiley--Interscience,  New York 1986.

\bibitem{Fel40}
W. Feller, {\it On the integro-differential equations of purely discontinuous
Markoff processes\/} Trans. Amer. Math. Soc., {\bf48} (1940), 488--815 and
Errata,  {\bf58} (1945) p. 474.

\bibitem{Fel57}
W. Feller, {\it On boundaries and lateral conditions for the Kolmogorov
differential equations\/}. Ann. Math. {\bf65} (1957), 527--570.

\bibitem{Fel59}
W. Feller, {\it The birth and death processes as diffusion processes\/}. J.
Math. Pures Appl. {\bf38} (1959), 301--345.

\bibitem{Feller-Book}
W. Feller, An introduction to probability theory and its applications. Vol. 1.
Wiley, 1970.

\bibitem{GT}
I. M. Gelfand, M. L.  Cetlin [Tsetlin],  {\it Finite-dimensional
representations of the group of unimodular matrices\/} (Russian). Doklady Akad.
Nauk SSSR (N.S.) {\bf71} (1950), no. 5, 825--828.


\bibitem{GnedinPitman}
A. Gnedin and J. Pitman, {\it Moment problems and boundaries of number
triangles\/}. Preprint, arXiv:0802.3410.

\bibitem{Gorin-FAA}
V. E. Gorin, {\it Disjointness of representations arising in the problem of
harmonic analysis on an infinite-dimensional unitary group\/}. Funktsional.
Anal. i Prilozhen. {\bf44} (2010), no. 2, 14--32 (Russian); translation in
Funct. Anal. Appl. {\bf44} (2010), no. 2, 92--105; arXiv:0805.2660.

\bibitem{GP13}
V. Gorin and G. Panova, {\it Asymptotics of symmetric polynomials with
applications to statistical mechanics and representation theory\/}. Preprint,
arXiv:1301.0634.

\bibitem{Ito06}
K. It\^o, {\it Essentials of stochastic processes\/}. Translated from the 1957
Japanese original. Translations of Mathematical Monographs, 231. American
Mathematical Society, Providence, RI, 2006.

\bibitem{ItoMcKean}

K. It\^o and H. P.  McKean, Jr., {\it Diffusion processes and their sample
paths\/}. Springer, 1965.

\bibitem{Karlin}
S. Karlin, {\it Total positivity I\/}. Stanford University Press, 1968.

\bibitem{KatoriTanemura-JSP09}
M. Katori and H. Tanemura, {\it Zeros of Airy function and relaxation
process\/}. J. Stat. Phys. {\bf136} (2009), 1177--1204; arXiv:0906.3666.

\bibitem{KatoriTanemura-CMP10}
M. Katori and H. Tanemura, {\it Non-equilibrium dynamics of Dyson's model with
an infinite number of particles\/}. Commun. Math. Phys. {\bf293} (2010),
469--497; arXiv:0812.4108

\bibitem{KatoriTanemura-MPRF11}
M. Katori and H. Tanemura, {\it Markov property of determinantal processes with
extended sine, Airy, and Bessel kernels\/}. Markov Processes and Relat. Fields
{\bf17} (2011), 541--580; arXiv:1106.4360.


\bibitem{Lig}
T. M. Liggett, {\it Continuous time Markov processes\/}. Graduate Texts in
Math. 113. Amer. Math. Soc., 2010.

\bibitem{McKean}
H. P. McKean, {\it Stochastic integrals\/}. Academic Press, 1969.

\bibitem{Mehta}
M. L. Mehta, {\it Random matrices\/}. Third edition. Academic Press, 2004.


\bibitem{Mey66}
P.-A. Meyer, {\it Probability and potentials\/}. Blaisdell, 1966.


\bibitem{OO98}
A. Okounkov and G. Olshanski, {\it Asymptotics of Jack polynomials as the
number of variables goes to infinity\/}. Intern. Math. Res. Notices {\bf1998}
(1998), no. 13, 641--682; arXiv:q-alg/9709011.

\bibitem{Ols-G&B}
G. Olshanski [Ol'shanskii], {\it Unitary representations of
infinite-dimensional pairs $(G,K)$ and the formalism of R. Howe\/}. In: {\it
Representation of Lie groups and related topics\/}. Advances in Contemp. Math.,
vol. 7 (A.~M.~Vershik and D.~P.~Zhelobenko, editors). Gordon and Breach, N.Y.,
London etc. 1990, 269--463.

\bibitem{Ols-AA89}
G. Olshanski, {\it  Unitary representations of $(G,K)$-pairs connected with the
infinite symmetric group $S(\infty)$\/}. Leningrad Math. J. 1, no. 4 (1990),
983--1014.

\bibitem{Ols-Semigroups}
G. Olshanski [Olshanskii], {\it On semigroups related to infinite-dimensional
groups\/}. In: {\it Topics in representation theory\/} (A.~A.~Kirillov, ed.).
Advances in Soviet Math., vol. 2. Amer. Math. Soc., Providence, R.I., 1991,
67--101.

\bibitem{Ols03}
G. Olshanski, {\it The problem of harmonic analysis on the infinite-dimensional
unitary group\/}. J. Funct. Anal. {\bf205} (2003), 464--524;
arXiv:math/0109193.

\bibitem{Osada-AnnProb13}
H. Osada, {\it Interacting Brownian motions in infinite dimensions with
logarithmic interaction potentials\/}. Ann. Prob. {\bf41} (2013), 1--49;
arXiv:0902.3561

\bibitem{Petrov}
L. Petrov, {\it The boundary of the Gelfand-Tsetlin graph: New proof of
Borodin-Olshanski's formula, and its q-analogue\/}. Moscow Math. J., to appear;
arXiv:1208.3443.


\bibitem{Pru63}
W. E. Pruitt, {\it Bilateral birth and death processes\/}. Trans. Amer. Math.
Soc. {\bf107} (1963), 508--525.

\bibitem{RogersWilliams}
L. C. G. Rogers and D. Williams, {\it Diffusions, Markov processes and
martingales\/}. Vols 1-2. Cambridge University Press, 2000.

\bibitem{Spo87}
H. Spohn, {\it Interacting Brownian particles: a study of Dyson's model\/}. In:
{\it Hydrodynamic Behavior and Interacting Particle Systems\/}, Papanicolaou,
G. (ed), IMA Volumes in Mathematics and its Applications, {\bf9}, Berlin:
Springer-Verlag, 1987, pp. 151--179.

\bibitem{VK81}
A. M. Vershik and S. V. Kerov, {\it Asymptotic theory of characters of the
symmetric group\/}. Funct. Anal. Appl. {\bf15} (1981), 246--255.

\bibitem{VK82}
A. M. Vershik and S. V. Kerov, {\it Characters and factor-representations of
the infinite unitary group\/}. Dokl. Akad. Nauk SSSR {\bf267} (1982), no. 2,
272--276 (Russian); English translation: Soviet Math. Dokl. {\bf26} (1982), no.
3, 570--574 (1983).

\bibitem{Voi76}
D. Voiculescu, {\it Repr\'esentations factorielles de type\/ {\rm II}$_1$ de
$U(\infty)$\/}. J. Math. Pures Appl. {\bf55} (1976), 1--20.


\bibitem{Wey39}
H. Weyl, {\it The classical groups. Their invariants and representations\/}.
Princeton Univ. Press, 1939; 1997 (fifth edition).

\bibitem{Winkler}
G. Winkler, {\it Choquet order and simplices. With applications in
probabilistic models\/}. Springer Lect. Notes Math. {\bf1145}, 1985.

\bibitem{Yan90}
Xiangqun Yang, {\it The construction theory of denumerable Markov processes\/}.
Hunan Science and Technology Publ. House, 1990 (Wiley series in probability and
mathematical statistics : Probability and mathematical statistics).

\bibitem{Zhe70}
D. P. Zhelobenko, {\it Compact Lie groups and their representations\/}, Nauka,
Moscow, 1970 (Russian); English translation: Transl. Math. Monographs 40, Amer.
Math. Soc., Providence, RI, 1973.




\end{thebibliography}
\end{document}